\documentclass[11pt]{article}
\usepackage{microtype}
\usepackage{fontspec}
\usepackage{lmodern}
\usepackage[british]{babel}

\usepackage{bm}
\usepackage{amssymb}
\usepackage{amsmath}
\usepackage{amsthm}
\usepackage{mathtools}
\usepackage{graphicx}
\usepackage[svgpath=figs/]{svg}
\usepackage{placeins}
\usepackage{multirow}
\usepackage[]{url}
\usepackage{float}
\usepackage[normalem]{ulem}
\usepackage{enumitem}

\usepackage[ruled,linesnumbered,noend]{algorithm2e}
\makeatletter
\newcommand{\RemoveAlgoNumber}{\renewcommand{\fnum@algocf}{\AlCapSty{\AlCapFnt\algorithmcfname}}}
\newcommand{\RevertAlgoNumber}{\algocf@resetfnum}
\makeatother
\usepackage{etoolbox}
\makeatletter
\renewcommand{\algocf@caption@ruled}{\AddToHook{cmd/MakeLinkTarget/before}
{\setlength\normalbaselineskip{\normalbaselineskip}}
\leavevmode\MakeLinkTarget[topalgocf]{algocf}\box\algocf@capbox\kern\interspacetitleruled\hrule
  width\algocf@ruledwidth height\algotitleheightrule depth0pt\kern\interspacealgoruled}%
\patchcmd\algocf@caption@algo
    {\hyper@refstepcounter{algocf}}
    {\hyper@refstepcounter{algocf}\xdef\@currentHref{top\@currentHref}}{}{\fail}
\makeatother

\usepackage{csquotes}
\usepackage[backend=biber,style=ieee,citestyle=numeric-comp,sorting=none,date=year,maxbibnames=99]{biblatex}
\usepackage{xpatch}
\xpatchbibdriver{online}
{\printtext[parens]{\usebibmacro{date}}}
{\iffieldundef{year}
    {}
    {\printtext[parens]{\usebibmacro{date}}}}
{}
{\typeout{There was an error patching biblatex-ieee (specifically, ieee.bbx's @online driver)}}
\AtEveryBibitem{%
  \ifentrytype{online}{}{% Remove url except for @online
    \clearfield{url}
  }
}

\usepackage{hyperref}
\hypersetup{
    colorlinks,
    citecolor={red},
    bookmarksdepth=3,
}

\usepackage[]{cleveref}
\crefname{equation}{Equation}{Equations}
\crefname{figure}{Figure}{Figures}
\crefname{section}{Section}{Sections}
\crefname{chapter}{Chapter}{Chapters}
\crefname{theorem}{Theorem}{Theorems}
\crefname{definition}{Definition}{Definitions}
\crefname{example}{Example}{Examples}
\crefname{algorithm}{Algorithm}{Algorithms}
\crefname{table}{Table}{Tables}
\crefname{appendix}{Appendix}{Appendices}

\usepackage{geometry}
\theoremstyle{plain}
\newtheorem{theorem}{Theorem}[section]

\newtheorem{lemma}[theorem]{Lemma}
\newtheorem{definition}[theorem]{Definition}

\theoremstyle{definition}

\DeclareMathOperator*{\argmin}{arg\,min}

\newcommand{\Vbeta}{\bm{g}^c}
\newcommand{\Vl}{\bm{l}}
\newcommand{\Vphi}{\bm{g}^f}
\newcommand\numberthis{\addtocounter{equation}{1}\tag{\theequation}}
\newcommand{\md}{\;\ifnum\currentgrouptype=16 \middle\fi|\;}

\newcommand{\Vb}{\bm{b}}

\newcommand{\Vg}{\bm{g}}

\newcommand{\Vo}{\bm{o}}
\newcommand{\Vp}{\bm{p}}

\newcommand{\Vu}{\bm{u}}
\newcommand{\Vv}{\bm{v}}
\newcommand{\Vw}{\bm{w}}
\newcommand{\Vx}{\bm{x}}
\newcommand{\Vy}{\bm{y}}
\newcommand{\Vz}{\bm{z}}
\newcommand{\MA}{\bm{A}}

\newcommand{\MI}{\bm{I}}

\newcommand{\ML}{\bm{L}}
\newcommand{\MM}{\bm{M}}
\newcommand{\MO}{\bm{O}}
\newcommand{\MP}{\bm{P}}

\newcommand{\MS}{\bm{S}}

\title{Accelerating MPGP-type Methods Through Preconditioning}
\author{Jakub Kružík \thanks{Department of Applied Mathematics, Faculty of Electrical Engineering and Computer Science, VSB-Technical University of Ostrava, Ostrava, Czech Republic}
\thanks{Institute of Geonics of the Czech Academy of Sciences, Ostrava, Czech Republic}
\thanks{Corresponding author: jakub.kruzik@vsb.cz}
\and David Horák$^{* \dag}$}

\begin{document}
\maketitle
\section*{Abstract}
This work investigates the acceleration of MPGP-type algorithms using preconditioning for the solution of quadratic programming problems.
The preconditioning needs to be done only on the free set so as not to change the constraints.
A variant of preconditioning restricted to the free set is the preconditioning in face.
The inner preconditioner in preconditioning in face needs to be recomputed or updated every time the free set changes.
Here, we investigate an approximate variant of preconditioning in face that computes the inner preconditioner only once.
We analyze the error of the approximate variant, give a sharp bound on the condition number of the preconditioned operator, and provide numerical experiments demonstrating that very large speedups can be achieved by the approximate variant.

\section{Introduction}
This work investigates the acceleration of the MPGP-type \cite{Dos-book-09} algorithms for the solution of quadratic programming (QP) problems
\begin{equation}
  \label{eq:qp}
  \argmin_{\Vx} \frac{1}{2}\Vx^{T}\MA\bm{x}-\Vx^{T}\Vb\quad\text{s.t.}\quad \Vx \in \Omega,
\end{equation}
where $\MA \in \mathbb{R}^{n\times n}$ is a symmetric positive semidefinite (SPS) matrix called the Hessian matrix, vector $\Vb \in \mathbb{R}^{n}$ is known as the right-hand side, and $\Omega$ is a set of constraints on the solution vector $\Vx \in \mathbb{R}^{n}$.
The minimized quadratic function $f(\Vx) = \frac{1}{2}\Vx^{T}\MA\bm{x}-\Vx^{T}\Vb$ is known as the cost function.
The MPGP-type algorithms employ projections onto the feasible set $\Omega$.
Therefore, the feasible set is typically assumed to be closed and convex so that the projection exists.
In our case, we will restrict ourselves to $\Omega$ consisting only of box constraints
\begin{equation*}
  \Omega = \left\{\Vx \in \mathbb{R}^n \mid \Vl \leq \Vx \leq \Vu \right\},
\end{equation*}
and the analysis will be done for positive definite Hessian.
%TODO l and u size is defined implicitly

The efficient solution of QP problems is important in a wide variety of fields, including economics, engineering, machine learning, and many others.
Concrete applications where MPGP-type algorithms were used include contact mechanics in fractured rock \cite{Stebel2024} with applications to modeling deep geological repositories of radioactive waste \cite{Claret2024}, machine learning for detecting wildfires from satellite images \cite{PechaPhD} or predicting compound bioactivity for the pharmaceutical industry \cite{Kruzik2018-svm}, and particle remapping for discrete element method modeling sea ice \cite{Turner2022}.

The restriction to the positive semidefinite Hessian allows us to use the conjugate gradient (CG) method for the unconstrained minimization part of the MPGP-type algorithm.
The CG method is very successful for the solution of large systems of linear equations, and many aspects of its convergence are well understood \cite{LiesenStrakos2012}.
Many problems solved by the CG method come from the discretization of PDEs with popular methods including finite elements/volumes, boundary elements, etc.
In these cases, refining the discretization worsens the conditioning of the resulting systems of linear equations, i.e., the Hessian, which results in the slowdown of the CG convergence.
The solution is to improve the spectrum of the Hessian using preconditioning.

Our aim is to modify MPGP-type algorithms to be able to use preconditioning in the underlying CG method while not changing the constraints.
This modification should not only lead to faster convergence in terms of the number of iterations but crucially in terms of time to solution.

The paper is divided as follows.
The following section briefly describes two MPGP-type algorithms - MPRGP and MPPCG.
In \cref{sec:prec}, we discuss the difficulty with preconditioning QP problems and show how preconditioning is implemented into the MPGP-type method.
\Cref{sec:facepre} describes preconditioning in face, and \cref{sec:approxpre} describes its approximate variant, including the error between the two variants and the conditioning of the preconditioned operator in a specific setting.
In \cref{sec:related}, we have enough prerequisites to describe related works.
Finally, we present numerical experiments in \cref{sec:results}.

\section{MPGP-type Algorithms}
\label{sec:mprgp}
QP problems with box constraints can be solved using the modified proportioning with reduced gradient projections (MPRGP) algorithm \cite{Dostal2005,Dos-book-09}.
The simplification to the feasible set with only one of the bound constraints is straightforward.
It is also possible to adapt the algorithm for various other constraints, such as elliptic and conical constraints \cite{Bouchala2014,PospisilPhD}.

As the name of the algorithm suggests, it utilizes gradient information for minimization, placing it among the first-order optimization methods.
While MPRGP does not directly work with active and free sets, the information about active and free sets is hidden in the gradient splitting, which is described later.
Consequently, MPRGP is considered an active set algorithm.
The algorithm was developed from the Polyak algorithm \cite{Polyak1969}.
%TODO note the diff with Polyak
A nice feature of the algorithm is that it has been proven to enjoy an R-linear rate of convergence given by the bound on the spectrum of the Hessian matrix \cite{Dos-book-09}.

In each iteration, MPRGP performs one of three types of steps: unconstrained minimization, expansion, or proportioning.
Since our Hessian $\MA$ is SPS, the unconstrained minimization is performed by a step of the CG method.
The active set is expanded by the \emph{expansion step}, which consists of a maximal feasible unconstrained minimization, in our case a partial CG step to the box, followed by a fixed step length gradient projection.
Finally, the \emph{proportioning step}, designed to reduce the active set, consists of a step of the steepest descent method in a proper direction.

To describe the algorithm, we first need to define the \emph{gradient splitting}.
Let $\bm{g}=\bm{A}\bm{x}-\bm{b}$ be the gradient of the cost function $f(\Vx)$ and let
\begin{equation*}
  \mathcal{A} =  \left\{ i\md \Vx_i=\Vl_i \quad\text{or}\quad \Vx_i=\Vu_i \right\}, \qquad
  \mathcal{F} =  \left\{ i\md \Vl_i <  \Vx_i < \Vu_i \right\}
\end{equation*}
be the \emph{active} and \emph{free set}, respectively.
Then the gradient splitting is defined component-wise for $i \in \{1, 2, \dots, n\}$ and is computed after each gradient evaluation.
The \emph{free gradient} $\Vg^f$ is defined as
\begin{flalign}
  \label{eq:gradSplit}
  \Vg_i^f = \begin{cases}
    0 \quad &\text{if}\quad i \in \mathcal{A},\\
    \Vg_i \quad &\text{if}\quad i \in \mathcal{F}.
    \end{cases}
\end{flalign}
%%TODO remove g^r?
%The \emph{reduced free gradient} $\Vg^r$ is defined as
%\begin{flalign*}
%   \Vg_i^r = \begin{cases} 0 \quad&\text{if}\quad i \in \mathcal{A},\\
%  \min\left(\frac{\Vx_i-\Vl_i}{\overline{\alpha}},\Vg_i\right) \quad &\text{if}\quad i \in \mathcal{F} \quad\text{and}\quad \Vg_i > 0,\\
%  \max\left(\frac{\Vx_i-\Vu_i}{\overline{\alpha}},\Vg_i\right)\quad &\text{if}\quad i \in \mathcal{F} \quad\text{and}\quad \Vg_i \le 0,
%  \end{cases}
%\end{flalign*}
%where $\overline{\alpha} \in (0,2||\bm{A}||^{-1}]$ is used as a priori chosen fixed step length in the expansion step~\cite{Dostal2005}.
%Typically, the parameter is given by
%\begin{equation*}
%\overline{\alpha} = \alpha_u ||\MA||^{-1},
%\end{equation*}
%where $\alpha_u \in (0,2]$ is a user-selected constant, and $||\MA||^{-1}$ is approximated by the power method.
%Effectively, $\Vphi$ is the gradient on the free set, and $\Vg^r$ is the free gradient that is reduced such that a steepest descent-type step in its opposite direction, that is $-\Vg^r$, with the step length $\overline{\alpha}$, does not leave the feasible set $\Omega$. A step in either of these directions can expand the active set but cannot reduce it.
A step in the direction $-\Vg^f$ may expand the active set but cannot reduce it.

The \emph{chopped gradient} $\Vg^c$ is defined as
\begin{flalign*}
  \Vg_i^c = \begin{cases} 0 \quad&\text{if}\quad i \in \mathcal{F},\\
  \min(\Vg_i,0) \quad&\text{if}\quad  \Vx_i = \Vl_i,\\
  \max(\Vg_i,0) \quad&\text{if}\quad \Vx_i = \Vu_i.
  \end{cases}
\end{flalign*}
A step in the direction $-\Vbeta$ may reduce the active set but cannot expand it.

The next ingredient is the projection onto the feasible set $\Omega$, which in the case of box constraints can be computed cheaply as
\begin{equation}
  \label{eq:proj}
\left[P_{\Omega}(\bm{x})\right]_i = \min\left\{\Vu_i,\max\left\{\Vl_i,\Vx_i\right\}\right\}, \quad i \in \{1, \dots, n\}.
\end{equation}

Finally, the \emph{projected gradient} is defined as $\bm{g}^{P}= \bm{g}^f + \bm{g}^c$.
The decrease in its norm serves as the natural stopping criterion for the algorithm since $\bm{g}^{P}= \bm{o}$ is equivalent to satisfying the Karush-Kuhn-Tucker conditions for our QP problem.

These are all the necessary ingredients to summarize MPRGP in Algorithm \ref{alg:mprgp}.

\begin{algorithm}[htbp]
    \SetStartEndCondition{ }{}{}%
    \DontPrintSemicolon
    \SetKwProg{Fn}{def}{\string:}{}
    \SetKwFunction{Range}{range}%%
    \SetKw{KwTo}{in}\SetKwFor{For}{for}{\string:}{}%
    \SetKwIF{If}{ElseIf}{Else}{if}{:}{elif}{else:}{}%
    \SetKwFor{While}{while}{:}{fintq}%
    \AlgoDontDisplayBlockMarkers\SetAlgoNoEnd\SetAlgoNoLine%
    %\scriptsize
    \KwIn{$\MA$, $\Vx_{0} \in \Omega$, $\Vb$, $\Gamma > 0$, $\overline{\alpha} \in (0,2||\MA||^{-1})$}
    $\Vg_0 = \MA \Vx_{0} - \Vb$, $\Vp_0 = \Vphi_0$, $k = 0$\;
    \While{$||\Vg^P_k||$ is not small}{
        \If{$||\Vbeta_k||^2 \le \Gamma^2 ||\Vg^f_k||^2$}{
         $\alpha^{feas}_k = \max \{\alpha\md \Vx_k - \alpha \Vp_k \in \Omega \}$\;
            $\alpha^{cg}_k = \Vg^T_k \Vp_k / \Vp^T_k \MA \Vp_k$\;
            \If{$\alpha^{cg}_k \le \alpha^{feas}_k$}{
              CG step - Algorithm \ref{alg:mpgpCg} \;
            }\Else{
              Expansion step - Algorithm \ref{alg:exp};
            }
        }\Else{
            Proportioning step - Algorithm \ref{alg:prop};
        }
        $k=k+1$\;
    }
    \KwOut{$\Vx_{k}$}

    \caption{MPRGP method}
    \label{alg:mprgp}
\end{algorithm}

\begin{algorithm}[htbp]
    \SetStartEndCondition{ }{}{}%
    \DontPrintSemicolon
    \SetKwProg{Fn}{def}{\string:}{}
    \SetKwFunction{Range}{range}%%
    \SetKw{KwTo}{in}\SetKwFor{For}{for}{\string:}{}%
    \SetKwIF{If}{ElseIf}{Else}{if}{:}{elif}{else:}{}%
    \SetKwFor{While}{while}{:}{fintq}%
    \AlgoDontDisplayBlockMarkers\SetAlgoNoEnd\SetAlgoNoLine%
    %\scriptsize
    $\Vx_{k+1} = \Vx_{k} -\alpha^{cg}_k\Vp_k$\;
    $\Vg_{k+1} = \Vg_k - \alpha^{cg}_k \MA \Vp_k$\;
    $\beta_k = \Vp^T_k \MA \Vphi_{k+1} / \Vp^T_k \MA \Vp_k$\;
    $\Vp_{k+1} = \Vphi_{k+1} -\beta_k \Vp_k$\;
    \caption{CG step}
    \label{alg:mpgpCg}
\end{algorithm}

\begin{algorithm}[htbp]
    \SetStartEndCondition{ }{}{}%
    \DontPrintSemicolon
    \SetKwProg{Fn}{def}{\string:}{}
    \SetKwFunction{Range}{range}%%
    \SetKw{KwTo}{in}\SetKwFor{For}{for}{\string:}{}%
    \SetKwIF{If}{ElseIf}{Else}{if}{:}{elif}{else:}{}%
    \SetKwFor{While}{while}{:}{fintq}%
    \AlgoDontDisplayBlockMarkers\SetAlgoNoEnd\SetAlgoNoLine%
    %\scriptsize
    $\Vx_{k+\frac{1}{2}} = \Vx_{k} - \alpha^{feas}_k \Vp_k$\;
    $\Vg_{k+\frac{1}{2}} = \Vg_k - \alpha^{feas}_k \MA \Vp_k$\;
    %$\Vx_{k+1} = P_{\Omega} (\Vx_{k+\frac{1}{2}} -\overline{\alpha} \Vphi(\Vx_{k+\frac{1}{2}}))$\;
    %$\Vx_{k+1} = P_{\Omega} (\Vx_{k+\frac{1}{2}} -\overline{\alpha} \Vphi) = \Vx_{k+\frac{1}{2}} -\overline{\alpha} \bm{g}^r$\;
    $\Vx_{k+1} = P_{\Omega} (\Vx_{k+\frac{1}{2}} -\overline{\alpha} \Vg_k^f)$\;
    $\Vg_{k+1} = \MA \Vx_{k+1} -\Vb$\;
    $\Vp_{k+1} = \Vphi_{k+1}$
    \caption{Expansion step}
    \label{alg:exp}
\end{algorithm}

\begin{algorithm}[htbp]
    \SetStartEndCondition{ }{}{}%
    \DontPrintSemicolon
    \SetKwProg{Fn}{def}{\string:}{}
    \SetKwFunction{Range}{range}%%
    \SetKw{KwTo}{in}\SetKwFor{For}{for}{\string:}{}%
    \SetKwIF{If}{ElseIf}{Else}{if}{:}{elif}{else:}{}%
    \SetKwFor{While}{while}{:}{fintq}%
    \AlgoDontDisplayBlockMarkers\SetAlgoNoEnd\SetAlgoNoLine%
    %\scriptsize
    $\alpha^{sd}_k = \Vg^T_k \Vbeta_k / (\Vbeta_k)^T \MA \Vbeta_k$\;
    $\alpha^{feas}_k = \max\{\alpha\md \Vx_k - \alpha \Vbeta_k \in \Omega\}$\;
    \If{$\alpha^{feas}_k < \alpha^{sd}_k$}{
      $\alpha^{sd}_k = \alpha^{feas}_k$\;
    }
    $\Vx_{k+1} = \Vx_{k} -\alpha^{sd}_k\Vbeta_k$\;
    $\Vg_{k+1} = \Vg_{k} - \alpha^{sd}_k \MA \Vbeta_k$\;
    $\Vp_{k+1} = \Vphi_{k+1}$\;
    \caption{Proportioning step}
    \label{alg:prop}
\end{algorithm}

%TODO check overline alpha intervals

The generalization to other constraints is in the way the maximal feasible step-length $\alpha_k^{feas}$ is computed, the definition and ease of computing the projection onto the feasible set $\Omega$, and potentially restricting the constant step-length to the first half of the interval, i.e., $\overline{\alpha} \in (0,||\MA||^{-1})$, when the set $\Omega$ is not subsymmetric \cite{Bouchala2014,Bouchala2012}.
The projections onto the feasible set $\Omega$ often have closed forms that are easily evaluated, as in our case given by \cref{eq:proj}.
Similarly, the computation of the maximal feasible step-length can also be very cheap.
In our case of the box-constraints, the closed formula is
\begin{equation*}
  \alpha_k^{feas} = \min \left\{\left(\Vx_i-\Vl_i\right)/\Vp_i : \Vp_i > 0, \min \left\{\left(\Vx_i-\Vu_i\right)/\Vp_i : \Vp_i < 0\right\} \right\},
\end{equation*}
where the minima are taken over all $i \in \{1,\dots,n\}$.
%TODO The convergence of the method for subsymmetric closed convex sets is established in.
%TODO infinite afeas, acg
%TODO proportioning afeas check

The MPRGP algorithm has Q-linear convergence in the cost function and, as a consequence, R-linear convergence in the norm of the error.
\begin{theorem}[R-linear convergence of MPRGP]
  Let $\Vx_k$ be generated by MRPGP with $\Vx_0 \in \Omega$, $\Gamma > 0$ and $\overline{\alpha} \in \left(0,2||\MA||^{-1}\right]$. Then
  \begin{equation*}
    f(\Vx_{k+1}) - f(\widehat{\Vx}) \le \eta\left( f(\Vx_k) - f(\widehat{\Vx}) \right)
    \quad \text{and} \quad ||\Vx_k - \widehat{\Vx}||_{\MA} \le 2 \eta^{k} \left(f(\Vx_0) - f(\widehat{\Vx})\right),
  \end{equation*}
  where $\widehat{\Vx}$ denotes the unique solution of \eqref{eq:qp} with $\Omega = \left\{\Vx \in \mathbb{R}^n \mid \Vl \leq \Vx \right\}$,
  \begin{displaymath}
    \eta = 1 - \frac{\widehat{\alpha}\lambda_{min}}{\vartheta \left(1+\widehat{\Gamma}^2\right)},
  \end{displaymath}
  \begin{displaymath}
    \widehat{\Gamma} = \max\{ \Gamma, \Gamma^{-1}\}, \qquad \vartheta = 2\max\{\overline{\alpha}||\MA||,1\}, \qquad \widehat{\alpha} = \min \{\overline{\alpha},2||\MA||^{-1} - \overline{\alpha}\}.
  \end{displaymath}
\end{theorem}
\begin{proof}
  See Theorem 5.14 in \cite{Dos-book-09}.
\end{proof}

The best convergence bound, which is achieved for $\Gamma = 1$ and $\overline{\alpha} = ||\MA||^{-1}$, reads
  \begin{displaymath}
    \eta^{opt} = 1 - \kappa \left(\MA \right)^{-1}/4,
  \end{displaymath}
where $\kappa \left(\MA \right)$ is the spectral condition number of $\MA$.

The MPRGP expansion step consists of the maximal feasible step in the direction of the CG direction that is followed by a fixed step-length gradient projection.
An improvement of the algorithm is to expand the active set using the full CG step that is projected, if needed, back to the feasible set.
The modified algorithm called modified proportioning with projected conjugate gradient (MPPCG) is obtained by replacing the expansion step \cref{alg:exp} with \cref{alg:projcg} in \cref{alg:mprgp}.
See \cite{Kruzik-projmprgp19,KruzikPhD} for more details and numerical comparison of MPPCG and MPRGP convergence speed.

\begin{algorithm}[htbp]
    \SetStartEndCondition{ }{}{}%
    \DontPrintSemicolon
    \SetKwProg{Fn}{def}{\string:}{}
    \SetKwFunction{Range}{range}%%
    \SetKw{KwTo}{in}\SetKwFor{For}{for}{\string:}{}%
    \SetKwIF{If}{ElseIf}{Else}{if}{:}{elif}{else:}{}%
    \SetKwFor{While}{while}{:}{fintq}%
    \AlgoDontDisplayBlockMarkers\SetAlgoNoEnd\SetAlgoNoLine%
    %\scriptsize
      $\Vx_{k+1} = P_{\Omega}(\Vx_{k} -\alpha^{cg}_k\Vp_k)$\;
      $\Vg_{k+1} = \MA \Vx_{k+1} -\Vb$\;
      $\Vp_{k+1} = \Vphi_{k+1}$\;
    \caption{Projected CG expansion step}
    \label{alg:projcg}
\end{algorithm}

\section{Preconditioned MPRGP and MPPCG Methods}
\label{sec:prec}
Preconditioning can significantly accelerate the CG method.
However, applying preconditioners to constrained QP problems is not straightforward.
Let us consider an SPD preconditioner matrix $\MM$ and the application of a preconditioner with this matrix as $\MM^{-1}$.
Using the split preconditioning to preserve the symmetry of the Hessian, the cost function is transformed into
\begin{equation*}
  f(\widehat{\Vx})_{\text{preconditioned}} = \frac{1}{2}\widehat{\bm{x}}^{T}\ML^{-1}\bm{A}\ML^{-T}\widehat{\Vx}-\widehat{\Vx}^{T}\ML^{-1}\bm{b},
\end{equation*}
where $\MM = \ML\ML^T$ and $\Vx = \ML^{-T}\widehat{\Vx}$.
Due to the variable change, the box constraints are transformed into general linear inequality constraints\footnote{Unless $\ML^{-T}$ is diagonal, e.g., when using a diagonal scaling preconditioner.}
\begin{equation*}
  \bm{l} \le \ML^{-T}\widehat{\Vx} \le \bm{u}.
\end{equation*}
QP problems with linear inequality constraints are typically much more difficult to solve.

Despite this, we will incorporate the preconditioning into MPGP-type methods disregarding the above disclaimer and only ensure that the constraints are not changed by the specific structure of the preconditioners, which are described in the following sections.
The preconditioning is incorporated into the MPGP-type methods in the same way as the preconditioning for the steepest descent and CG methods is incorporated; see e.g. \cite{golub13,LiesenStrakos2012}.
Denoting $\MM^{-1}$ as the preconditioner action, then the preconditioned MPRGP algorithm can be found in \cref{alg:mprgpPre}.
The preconditioned MPPCG method is obtained by replacing the preconditioned expansion step (\cref{alg:expPre}) with the preconditioned projected CG step (\cref{alg:projcgPre}) in \cref{alg:mprgpPre}.

\begin{algorithm}[htbp]
    \SetStartEndCondition{ }{}{}%
    \DontPrintSemicolon
    \SetKwProg{Fn}{def}{\string:}{}
    \SetKwFunction{Range}{range}%%
    \SetKw{KwTo}{in}\SetKwFor{For}{for}{\string:}{}%
    \SetKwIF{If}{ElseIf}{Else}{if}{:}{elif}{else:}{}%
    \SetKwFor{While}{while}{:}{fintq}%
    \AlgoDontDisplayBlockMarkers\SetAlgoNoEnd\SetAlgoNoLine%
    %\scriptsize
    \KwIn{$\MA$, $\MM^{-1}$, $\Vx_{0} \in \Omega$, $\Vb$, $\Gamma > 0$, $\overline{\alpha} \in (0,2||\MA||^{-1})$}
    $\Vg_0 = \MA \Vx_{0} - \Vb$, $\Vz_0 = \MM^{-1}\Vphi_0$, $\Vp_0 = \Vz_0$, $k = 0$\;
    \While{$||\Vg^P_k||$ is not small}{
        \If{$||\Vbeta_k||^2 \le \Gamma^2 ||\Vg^f_k||^2$}{
         $\alpha^{feas}_k = \max \{\alpha \md \Vx_k - \alpha \Vp_k \in \Omega \}$\;
            $\alpha^{cg}_k = \Vg^T_k \Vz_k / \Vp^T_k \MA \Vp_k$\;
            \If{$\alpha^{cg}_k \le \alpha^{feas}_k$}{
              Preconditioned CG step - Algorithm \ref{alg:cgPre} \;
            }\Else{
              Preconditioned expansion step - Algorithm \ref{alg:expPre};
            }
        }\Else{
            Preconditioned proportioning step - Algorithm \ref{alg:propPre};
        }
        $k=k+1$\;
    }
    \KwOut{$\Vx_{k}$}

    \caption{Preconditioned MPRGP}
    \label{alg:mprgpPre}
\end{algorithm}

\begin{algorithm}[htbp]
    \SetStartEndCondition{ }{}{}%
    \DontPrintSemicolon
    \SetKwProg{Fn}{def}{\string:}{}
    \SetKwFunction{Range}{range}%%
    \SetKw{KwTo}{in}\SetKwFor{For}{for}{\string:}{}%
    \SetKwIF{If}{ElseIf}{Else}{if}{:}{elif}{else:}{}%
    \SetKwFor{While}{while}{:}{fintq}%
    \AlgoDontDisplayBlockMarkers\SetAlgoNoEnd\SetAlgoNoLine%
    %\scriptsize
    $\Vx_{k+1} = \Vx_{k} -\alpha^{cg}_k\Vp_k$\;
    $\Vg_{k+1} = \Vg_k - \alpha^{cg}_k \MA \Vp_k$\;
    $\Vz_{k+1} = \MM^{-1} \Vphi_{k+1}$\;
    $\beta_k = \Vp^T_k \MA \Vz_{k+1} / \Vp^T_k \MA \Vp_k$\;
    $\Vp_{k+1} = \Vz_{k+1} -\beta_k \Vp_{k}$\;
    \caption{Preconditioned CG step}
    \label{alg:cgPre}
\end{algorithm}

\begin{algorithm}[htbp]
    \SetStartEndCondition{ }{}{}%
    \DontPrintSemicolon
    \SetKwProg{Fn}{def}{\string:}{}
    \SetKwFunction{Range}{range}%%
    \SetKw{KwTo}{in}\SetKwFor{For}{for}{\string:}{}%
    \SetKwIF{If}{ElseIf}{Else}{if}{:}{elif}{else:}{}%
    \SetKwFor{While}{while}{:}{fintq}%
    \AlgoDontDisplayBlockMarkers\SetAlgoNoEnd\SetAlgoNoLine%
    %\scriptsize
    $\Vx_{k+\frac{1}{2}} = \Vx_{k} - \alpha^{feas}_k \Vp_k$\;
    $\Vg_{k+\frac{1}{2}} = \Vg_k - \alpha^{feas}_k \MA \Vp_k$\;
    %$\Vx_{k+1} = P_{\Omega} (\Vx_{k+\frac{1}{2}} -\overline{\alpha} \Vphi(\Vx_{k+\frac{1}{2}}))$\;
    %$\Vx_{k+1} = P_{\Omega} (\Vx_{k+\frac{1}{2}} -\overline{\alpha} \Vphi) = \Vx_{k+\frac{1}{2}} -\overline{\alpha} \bm{g}^r$\;
    $\Vx_{k+1} = P_{\Omega} (\Vx_{k+\frac{1}{2}} -\overline{\alpha} \Vg_k^f)$\;
    $\Vg_{k+1} = \MA \Vx_{k+1} -\Vb$\;
    $\Vz_{k+1} = \MM^{-1} \Vphi_{k+1}$\;
    $\Vp_{k+1} = \Vz_{k+1}$\;
    \caption{Preconditioned expansion step}
    \label{alg:expPre}
\end{algorithm}

\begin{algorithm}[htbp]
    \SetStartEndCondition{ }{}{}%
    \DontPrintSemicolon
    \SetKwProg{Fn}{def}{\string:}{}
    \SetKwFunction{Range}{range}%%
    \SetKw{KwTo}{in}\SetKwFor{For}{for}{\string:}{}%
    \SetKwIF{If}{ElseIf}{Else}{if}{:}{elif}{else:}{}%
    \SetKwFor{While}{while}{:}{fintq}%
    \AlgoDontDisplayBlockMarkers\SetAlgoNoEnd\SetAlgoNoLine%
    %\scriptsize
    $\alpha^{sd}_k = \Vg^T_k \Vbeta_k / (\Vbeta_k)^T \MA \Vbeta_k$\;
    $\alpha^{feas}_k = \max\{\alpha \md \Vx_k - \alpha \Vbeta_k \in \Omega\}$\;
    \If{$\alpha^{feas}_k < \alpha^{sd}_k$}{
      $\alpha^{sd}_k = \alpha^{feas}_k$\;
    }
    $\Vx_{k+1} = \Vx_{k} -\alpha^{sd}_k\Vbeta_k$\;
    $\Vg_{k+1} = \Vg_{k} - \alpha^{sd}_k \MA \Vbeta_k$\;
    $\Vz_{k+1} = \MM^{-1} \Vphi_{k+1}$\;
    $\Vp_{k+1} = \Vz_{k+1}$\;
    \caption{Preconditioned proportioning step}
    \label{alg:propPre}
\end{algorithm}

\begin{algorithm}[htbp]
    \SetStartEndCondition{ }{}{}%
    \DontPrintSemicolon
    \SetKwProg{Fn}{def}{\string:}{}
    \SetKwFunction{Range}{range}%%
    \SetKw{KwTo}{in}\SetKwFor{For}{for}{\string:}{}%
    \SetKwIF{If}{ElseIf}{Else}{if}{:}{elif}{else:}{}%
    \SetKwFor{While}{while}{:}{fintq}%
    \AlgoDontDisplayBlockMarkers\SetAlgoNoEnd\SetAlgoNoLine%
    %\scriptsize
    $\Vx_{k+1} = P_{\Omega}(\Vx_{k} -\alpha^{cg}_k\Vp_k)$\;
    $\Vg_{k+1} = \MA \Vx_{k+1} -\Vb$\;
    $\Vz_{k+1} = \MM^{-1} \Vphi_{k+1}$\;
    $\Vp_{k+1} = \Vz_{k+1}$\;
    \caption{Preconditioned projected CG step}
    \label{alg:projcgPre}
\end{algorithm}

\section{Preconditioning in Face}
\label{sec:facepre}
Preconditioning in face was introduced in \cite{Oleary1980} for the Polyak algorithm, and its use is described for the MPRGP algorithm in \cite{Dos-book-09}.

The idea is to apply the preconditioning only on the free set.
%TODO where the constraints are not binding...!!!!!!!!!
In order to achieve this, we split the preconditioner matrix according to the free set and the active set
%TODO explain overline M
\begin{equation}
  \label{eq:split}
  \overline{\MM} = \begin{pmatrix}
   \MM_{\mathcal{F}\mathcal{F}} & \MM_{\mathcal{F}\mathcal{A}} \\
   \MM_{\mathcal{A}\mathcal{F}} & \MM_{\mathcal{A}\mathcal{A}}
 \end{pmatrix}.
\end{equation}
Then only the free gradient is preconditioned by a preconditioner computed solely on the free set
\begin{equation}
  \label{eq:prec}
  \Vz = \begin{pmatrix}
    \Vz^f_\mathcal{F}\\
    \Vo
  \end{pmatrix} =
  \MM^{-1}
 \begin{pmatrix}
 \Vg^f_\mathcal{F} \\ \Vo
 \end{pmatrix} \coloneq
 \begin{pmatrix}
   \MM_{\mathcal{F}\mathcal{F}}^{-1} & \Vo \\
   \Vo & \Vo
 \end{pmatrix}
 \begin{pmatrix}
 \Vg^f_\mathcal{F} \\ \Vo
 \end{pmatrix},
\end{equation}
where $\MM^{-1}$ is the application of the preconditioning in face, while $\MM_{\mathcal{F}\mathcal{F}}^{-1}$ is an application of some standard preconditioner like incomplete Cholesky.
We call $\MM_{\mathcal{F}\mathcal{F}}^{-1}$ the application of the inner preconditioner.
Notice that the preconditioning in face gives something like a preconditioned free gradient $\Vz^f$.
We note that the vectors are usually not reordered in actual implementations.

Obviously, the major drawback is that the preconditioner must be recomputed or at least updated every time the free set changes.
One way to avoid recomputing the preconditioner is to restrict the preconditioner not to the current free set, but to the set that will never be active.
Then the preconditioner needs to be computed only once.
For example, if only a part of the solution vector is constrained, the preconditioner can be computed and applied only to the unconstrained part.
Such problems arise in, e.g., contact problems.
For example, let us take the case of the 3D cube with a contact interface on only one of its sides, which is described in more detail in later \cref{sec:results}.
If the cube is discretized with $n \times n \times n$ unknowns, only $n^{2}$ unknowns, i.e., at most $1/n$ of all unknowns, can become active.
This allows us to apply preconditioning to the remaining $\left(n-1\right)n^2$ unknowns.

In the following section, we develop an alternative preconditioning method that avoids the need to recompute the preconditioner without prior knowledge of the set that will never be active.

\section{Approximate Preconditioning in Face}
\label{sec:approxpre}
To avoid the need to recompute the preconditioner, it is possible to apply the full preconditioner, which is denoted $\overline{\MM}^{-1}$, computed for the entire preconditioning matrix $\overline{\MM}$, and then zero out the active set components
\begin{flalign*}
  \Vz = \begin{pmatrix}
  \widetilde{\Vz}^f_\mathcal{F}\\
    \Vo
  \end{pmatrix} =
  \MM^{-1}
 \begin{pmatrix}
 \Vg^f_\mathcal{F} \\ \Vo
\end{pmatrix}
  &\coloneq \text{gradientSplit}_{Free}( \overline{\MM}^{-1}
 \begin{pmatrix}
 \Vg^f_\mathcal{F} \\ \Vo
\end{pmatrix}),
\end{flalign*}
where the function $\text{gradientSplit}_{Free}()$ zeros out the active set components of a given vector in the same way as computing the free gradient in \cref{eq:gradSplit}.
The operator $\MM^{-1}$ is again the application of the preconditioner, which we call the \emph{approximate preconditioning in face} because it tries to approximate the preconditioned free gradient $\Vz^f$ from \cref{eq:prec}.
The operator $\overline{\MM}^{-1}$ is the application of some standard preconditioner, disregarding any information about the free set. Again we call $\overline{\MM}^{-1}$ the inner preconditioner.
Assuming that $\Vg^{f}$ is split consistently with \eqref{eq:split}, we can define the application of the approximate preconditioner in face equivalently as
\begin{equation*}
  \MM^{-1} \Vg^{f} \coloneq \MP_{\mathcal{F}} \overline{\MM}^{-1} \MP_{\mathcal{F}} \Vg^{f}, \quad\text{where}\quad \MP_{\mathcal{F}} = \begin{pmatrix} \MI_{\mathcal{FF}} & \MO\\ \MO & \MO \end{pmatrix}.
\end{equation*}

Let us assume for the rest of the section that $\overline{\MM} = \MA$ and the application of the inner preconditioner $\overline{\MM}^{-1}$ is the actual inverse.
Note that this means that any matrix inverse notation for the rest of this section also represents the inverse of a matrix and not an application of some preconditioner.
These assumptions will allow us to establish some interesting properties of the approximate preconditioning in face.

With the above assumptions, the approximate preconditioner corresponds to the preconditioning by the Schur complement eliminating the active set variables.
\begin{theorem}[Schur complement relation]
  \label{theo:schur}
The application of the approximate preconditioning in face corresponds to the multiplication of the free gradient restricted to the free set by the inverse of the Schur complement
\begin{equation*}
\bm{S}^{-1} = (\MM_{\mathcal{F}\mathcal{F}} - \MM_{\mathcal{F}\mathcal{A}} \MM_{\mathcal{A}\mathcal{A}}^{-1} \MM_{\mathcal{A}\mathcal{F}})^{-1},
\end{equation*}
and padding the resulting vector by zeros on the active set.

Moreover, the preconditioned operator $\MM^{-1}\MA$ is
\begin{equation}
  \label{eq:precondop}
  \MM^{-1}\MA = \begin{pmatrix}
  \bm{S}^{-1}\MA_{\mathcal{F}\mathcal{F}} & \MS^{-1}\MA_{\mathcal{FA}}\\
  \MO & \MO
\end{pmatrix}.
\end{equation}

%Therefore, the effective condition number of the preconditoned operator is
%\begin{equation*}
%  \kappa_{eff}\left( \MM^{-1}\MA \right) =
%\end{equation*}
\end{theorem}
\begin{proof}
Using the block inversion formula, we have
\begin{equation*}
\overline{\MM}^{-1} = \begin{pmatrix}
  \MS^{-1} & -\MS^{-1}\MA_{\mathcal{FA}}\MA_{\mathcal{AA}}^{-1}\\
  -\MA_{\mathcal{AA}}^{-1}\MA_{\mathcal{AF}}\MS^{-1} & \MA_{\mathcal{AA}}^{-1} + \MA_{\mathcal{AA}}^{-1}\MA_{\mathcal{AF}}\MS^{-1}\MA_{\mathcal{FA}}\MA_{\mathcal{AA}}^{-1}
\end{pmatrix},
\end{equation*}
and substituting into the preconditioner definition yields
\begin{flalign*}
  \text{gradientSplit}_{Free}( \overline{\MM}^{-1}
 \begin{pmatrix}
 \Vg^f_\mathcal{F} \\ \Vo
\end{pmatrix}) =
 \begin{pmatrix}
   (\MM_{\mathcal{F}\mathcal{F}} - \MM_{\mathcal{F}\mathcal{A}} \MM_{\mathcal{A}\mathcal{A}}^{-1} \MM_{\mathcal{A}\mathcal{F}})^{-1} \Vg^f_\mathcal{F} \\ \Vo
 \end{pmatrix} =
 \begin{pmatrix}
   \bm{S}^{-1} \Vg^f_\mathcal{F} \\ \Vo
 \end{pmatrix}.
\end{flalign*}

When applying the preconditioner to $\MA$, assuming $\MA$ is partitioned as in \eqref{eq:split}, we have using the block inversion formula above
\begin{equation*}
\MM^{-1}\MA = \MP_{\mathcal{F}} \overline{\MM}^{-1} \MP_{\mathcal{F}}\MA =
\begin{pmatrix} \MS^{-1} & \MO\\ \MO & \MO \end{pmatrix} \MA =
\begin{pmatrix} \MS^{-1}\MA_{\mathcal{FF}} & \MS^{-1}\MA_{\mathcal{FA}}\\ \MO & \MO \end{pmatrix}.
\end{equation*}
\end{proof}

With the above characterization of our preconditioned operator, we can give an error term of the approximate variant preconditioner as compared to the standard preconditioning in face.
\begin{theorem}[The error of the approximate preconditioning in face]
  \label{theo:eigvals}
  The (TODO additive?) error of the approximate preconditioning in face compared to the standard preconditioning in face is given by
\begin{equation*}
  \MM_{\mathcal{F}\mathcal{F}}^{-1} \MM_{\mathcal{F}\mathcal{A}} ( \MM_{\mathcal{A}\mathcal{A}} - \MM_{\mathcal{A}\mathcal{F}} \MM_{\mathcal{F}\mathcal{F}}^{-1} \MM_{\mathcal{F}\mathcal{A}})^{-1} \MM_{\mathcal{A}\mathcal{F}}.
\end{equation*}

Moreover, the non-zero eigenvalues of the preconditioned operator $\MM^{-1}\MA$ correspond to the eigenvalues of $\MS^{-1}\MA_{\mathcal{FF}}$, which are
\begin{equation*}
  1 = \lambda_1 = \dots = \lambda_{m-r} \le \dots \le \lambda_m,
\end{equation*}
where $r=\text{rank}(\MM_{\mathcal{AF}})$ and $m$ is the size of the free set.
\end{theorem}
\begin{proof}
  From \eqref{eq:precondop} the only possibly non-zero eigenvalues of $\MM^{-1}\MA$ are eigenvalues of $\MS^{-1}\MA_{\mathcal{FF}}$.
  Using the alternative matrix block inversion formula, we have
\begin{flalign*}
  \bm{S}^{-1} \MA_{\mathcal{FF}} &=
   (\MM_{\mathcal{F}\mathcal{F}}^{-1} + \MM_{\mathcal{F}\mathcal{F}}^{-1} \MM_{\mathcal{F}\mathcal{A}} ( \MM_{\mathcal{A}\mathcal{A}} - \MM_{\mathcal{A}\mathcal{F}} \MM_{\mathcal{F}\mathcal{F}}^{-1} \MM_{\mathcal{F}\mathcal{A}})^{-1} \MM_{\mathcal{A}\mathcal{F}} \MM_{\mathcal{F}\mathcal{F}}^{-1}) \MA_{\mathcal{FF}}
 \\&=
   (\MI + \MM_{\mathcal{F}\mathcal{F}}^{-1} \MM_{\mathcal{F}\mathcal{A}} ( \MM_{\mathcal{A}\mathcal{A}} - \MM_{\mathcal{A}\mathcal{F}} \MM_{\mathcal{F}\mathcal{F}}^{-1} \MM_{\mathcal{F}\mathcal{A}})^{-1} \MM_{\mathcal{A}\mathcal{F}}) \MM_{\mathcal{F}\mathcal{F}}^{-1} \MA_{\mathcal{FF}}
 \\&=
   \MI + \MM_{\mathcal{F}\mathcal{F}}^{-1} \MM_{\mathcal{F}\mathcal{A}} ( \MM_{\mathcal{A}\mathcal{A}} - \MM_{\mathcal{A}\mathcal{F}} \MM_{\mathcal{F}\mathcal{F}}^{-1} \MM_{\mathcal{F}\mathcal{A}})^{-1} \MM_{\mathcal{A}\mathcal{F}}.
  \numberthis \label{eq:precondexpand}
\end{flalign*}
Since the preconditioned operator for the preconditioning in face would be equal to $\MI$, we have that the term
\begin{equation*}
  \MM_{\mathcal{F}\mathcal{F}}^{-1} \MM_{\mathcal{F}\mathcal{A}} ( \MM_{\mathcal{A}\mathcal{A}} -                       \MM_{\mathcal{A}\mathcal{F}} \MM_{\mathcal{F}\mathcal{F}}^{-1} \MM_{\mathcal{F}\mathcal{A}})^{-1} \MM_{\mathcal{A}\mathcal{F}}
\end{equation*}
is the error of the approximate preconditioning in face compared to the standard preconditioning in face.

Moreover, the eigenvalues of $\MS^{-1}\MA_{\mathcal{FF}}$ are by \eqref{eq:precondexpand} non-zero and equal to
\begin{equation*}
   1 = \lambda_1 = \dots = \lambda_{m-r} \le \dots \le \lambda_m,
\end{equation*}
where $r=\text{rank}(\MM_{\mathcal{AF}})$ and $m$ is the size of the free set.
\end{proof}
As can be seen from the previous theorem, the difference between the two preconditioners is only in at most $\text{rank}(\MM_{\mathcal{AF}})$ eigenvalues.

\begin{figure}[htbp]
\includegraphics[width=1.0\textwidth]{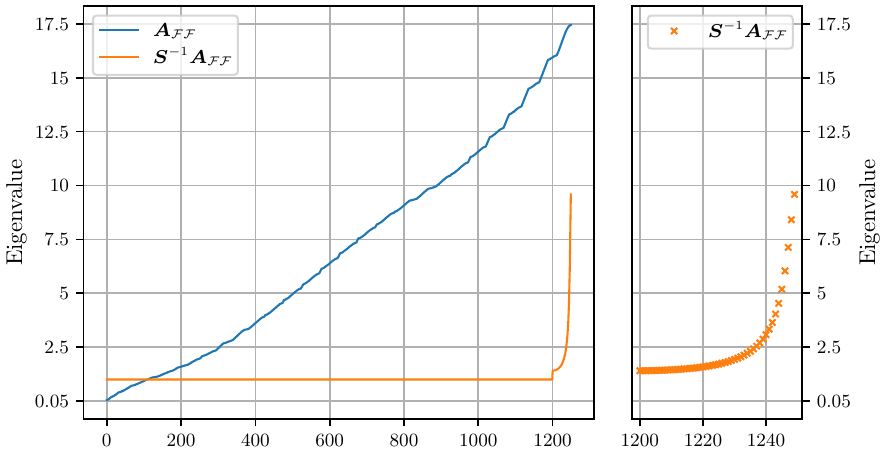}
\caption[Eigenvalues of the journal bearing problem]{Eigenvalues of the journal bearing problem with 50x50 grid points (2,500 DOFs) preconditioned by the inverse matrix at iteration 0 (the free set size is 1,250, and the rank of the off-diagonal block is 50).}
\label{fig:eigsPrecond}
\end{figure}

To illustrate \cref{theo:eigvals}, we plotted in \cref{fig:eigsPrecond} the eigenvalues for the journal bearing problem, which is described later in \cref{sec:results}, in the first iteration with the zero initial guess.
The difference between the preconditioning in face and its approximate variant is precisely in the last 50 eigenvalues, since $\text{rank}(\MM_{\mathcal{AF}}) = 50$.
We note that those last 50 eigenvalues are spaced throughout the interval starting at 1 and ending with some maximal eigenvalue.

\begin{figure}[htbp]
\centering
\includegraphics[width=1.\columnwidth]{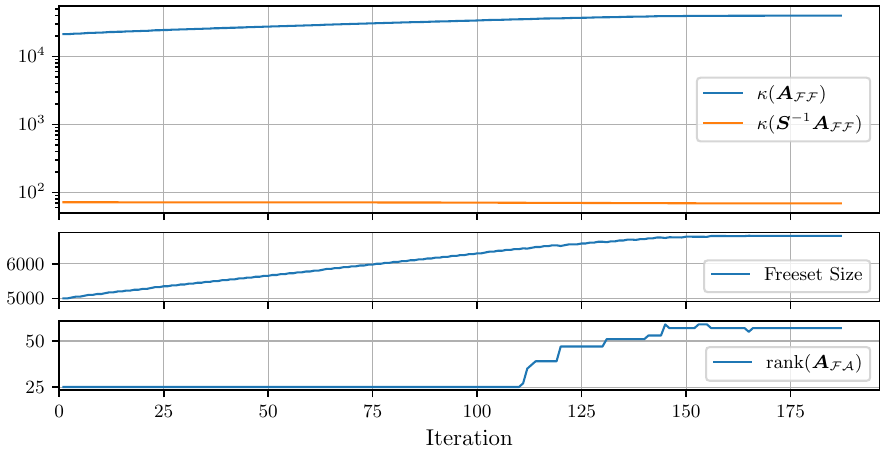}
\caption[Preconditioned journal bearing problem]{Condition number, free set size, and rank of the off-diagonal block for the journal bearing problem with 400x25 grid points (10,000 DOFs) preconditioned by the inverse matrix.}
\label{fig:precondCondJbearing}
\end{figure}

To see how the condition number and the rank of the off-diagonal matrix $\MM_{\mathcal{AF}}$ change throughout the iterative process, we plotted these quantities together with the free set size for the journal bearing problem with a different discretization in \cref{fig:precondCondJbearing}.
We can see that the condition number of the preconditioned operator remains essentially constant and that it was always significantly lower than the condition number of the unpreconditioned operator.
The off-diagonal matrix rank grew moderately from $25$ to the maximum of $59$ for the journal bearing problem, which represented only a tiny fraction of the free set size where the preconditioning is applied.

%\begin{figure}[htbp]
%\centering
%\includegraphics[width=1.\columnwidth]{graphs/ex1_cond.pdf}
%\caption[Preconditioned 1D Poisson's problem]{Condition number, free set size, and rank of the off-diagonal block for \emph{ex1} with 1,000 DOFs preconditioned by the inverse matrix.}
%\label{fig:precondCondEx1}
%\end{figure}

In the following, we bound the condition number of the preconditioned operator under our assumptions of $\overline{\MM} = \MA$ and the application of the inner preconditioner $\overline{\MM}^{-1}$ is the actual inverse. First, we need the Cauchy--Bunyakowski--Schwarz (CBS) constant appearing in the CBS inequality.
\begin{definition}[The CBS constant (Section 9.1 in \cite{Axelsson96})]
  \label{theo:cbs}
  The CBS constant $\gamma$ is defined as
  \begin{equation*}
    \gamma = \sup_{\Vw_1 \in W_1, \Vw_2 \in W_2} \frac{\Vw_1 \MA \Vw_2}{\left(\Vw_1^{T}\MA\Vw_1 \Vw_2^{T}\MA\Vw_2\right)^{\frac{1}{2}}},
  \end{equation*}
  where
  \begin{equation*}
    W_1 = \left\{ \Vv = \begin{pmatrix} \Vv_\mathcal{F}\\ \Vo  \end{pmatrix}, \Vv \in \mathbb{R}^{n} \right\}, \quad W_2 = \left\{ \Vv = \begin{pmatrix} \Vo \\ \Vv_\mathcal{A}  \end{pmatrix}, \Vv \in \mathbb{R}^{n} \right\}.
  \end{equation*}
\end{definition}

The CBS constant can be used to estimate the conditioning of our preconditioned operator $\MS^{-1}\MA_{\mathcal{FF}}$.
  \begin{lemma}
    \label{theo:schurcond}
    For all $\Vv_{\mathcal{F}} \neq \Vo$
    \begin{equation*}
      1 \leq \frac{\Vv_{\mathcal{F}}^{T} \MA_{\mathcal{FF}} \Vv_{\mathcal{F}}}{\Vv_{\mathcal{F}}^{T} \MS \Vv_{\mathcal{F}}} \leq \frac{1}{1-\gamma^{2}}.
    \end{equation*}
    The right-hand inequality is sharp and the left-hand inequality is sharp if $\mathcal{N}(\MA_{\mathcal{AF}}$) is non-trivial.
  \end{lemma}
  \begin{proof}
    Follow the proof of Lemma 9.2 (b) in \cite{Axelsson96} and take the reciprocal of the result. (Axelsson does the proof for the active set in our notation - simply permuting the diagonal blocks gives the result on the free set).
  \end{proof}

The Wielandt inequality \cite{Wielandt1953} will give us the upper bound on the CBS constant $\gamma$.
\begin{theorem}[Wielandt inequality]
  \label{theo:wielandt}
  Let $\lambda_{min}$ and $\lambda_{max}$ be the smallest and the largest eigenvalue of $\MA$, respectively. Then for all $\Vx, \Vy \in \mathbb{R}^{n}$ such that $\Vx^{T}\Vy = 0$ it holds that
  \begin{equation*}
    |\Vx^T \MA \Vy|^{2} \leq \left( \frac{\lambda_{min}-\lambda_{max}}{\lambda_{min}+\lambda_{max}} \right)^{2} \left( \Vx^{T}\MA\Vx \right) \left( \Vy^{T} \MA \Vy \right).
  \end{equation*}
\end{theorem}
\begin{proof}
 See Section 7.4.12 of \cite{Horn2012}.
\end{proof}

Now we are ready to prove the condition number of the preconditioned operator based on the condition number of the Hessian $\MA$.
\begin{theorem}{Conditioning of the preconditioned operator}
  \label{theo:precond}
  Let $\overline{\MM} = \MA$ and the inner preconditioner is the full inverse of $\overline{\MM}$. Then
  %TODO K \le 1/(1-gamma^2) \le ...
    \begin{equation*}
      \kappa\left(\MM^{-1}\MA\right) \leq \frac{\left( \kappa\left( \MA \right) +1  \right)^{2}}{4\kappa(\MA)}.
    \end{equation*}
\end{theorem}
\begin{proof}
  Since $\MM^{-1} \MA = \MS^{-1}\MA_{\mathcal{FF}}$, we have from Lemma \ref{theo:schurcond}
  \begin{equation}
    \label{eq:condcbs}
    \kappa\left(\MM^{-1}\MA\right) \leq \frac{1}{1-\gamma^{2}}.
  \end{equation}
  Combining Definition \ref{theo:cbs} and Theorem \ref{theo:wielandt} by using $\Vx \in W_1$ and $\Vy \in W_2$, we get the estimate of the square of the CBS constant as
  \begin{equation*}
    \gamma^{2} \leq \left( \frac{\lambda_{min}-\lambda_{max}}{\lambda_{min}+\lambda_{max}} \right)^{2} = \left(\frac{\kappa(A)-1}{\kappa(A)+1}\right)^{2}.
  \end{equation*}
  Using the $\gamma^{2}$ estimate in \eqref{eq:condcbs} gives
  \begin{equation*}
    \kappa\left(\MM^{-1}\MA\right) \leq \frac{1}{1-\left(\frac{\kappa(\MA)-1}{\kappa(\MA)+1}\right)^{2}} = \frac{1}{\frac{\left(\kappa(\MA)+1\right)^{2} - \left(\kappa(\MA)-1\right)^{2}}{\left(\kappa(\MA)+1\right)^{2}}} = \frac{1}{\frac{4\kappa(\MA)}{\left(\kappa(\MA)+1\right)^{2}}} = \frac{\left( \kappa\left( \MA \right) +1  \right)^{2}}{4\kappa(\MA)}
  \end{equation*}
\end{proof}

Since by equation \eqref{eq:condcbs} our condition number estimate depends on the CBS constant $\gamma$, we know that $\gamma \in \left[0, 1\right)$ (for SPD matrix $\MA$) \cite{Axelsson96} and that the constant describes the "strength" of coupling between the diagonal blocks or equivalently between active and free sets in our case.
For fully decoupled blocks, we have $\gamma = 0$ so that by \eqref{eq:condcbs} we have $\kappa\left( \MM^{-1}\MA \right) = 1$.
For instance, let
\begin{equation*}
  \MA = \begin{pmatrix} 2 & 0 \\ 0 & 1 \end{pmatrix} \qquad \text{and} \qquad \mathcal{F} = \{1, 2\},
\end{equation*}
then indeed $\kappa\left( \MM^{-1}\MA \right) = 1$.
However, our estimate gives $\kappa\left( \MM^{-1}\MA \right) \leq 9/8$, and replacing the $2$ in $\MA$ with a very large number like $1$ million will severely overestimate the preconditioned operator condition number.
Let us construct an example with strong coupling, i.e., $\gamma$ is close to $1$. Let $n \geq 3$ be odd, the Hessian be the Laplace $n \times n$ tridiagonal matrix
\begin{equation*}
  \MA = \begin{pmatrix}
    2 & -1 & 0 & 0 & \dots & 0 \\
    -1 & 2 & -1 & 0 & \dots & 0 \\
     \ddots & \ddots & \ddots & \ddots & \ddots & \ddots
   \end{pmatrix},
\end{equation*}
and the free set consist of odd numbers between $1$ and $n$.
Then \eqref{eq:condcbs} is sharp because $\mathcal{N}(\MA_{\mathcal{AF}}$) is non-trivial.
%The CBS constant $\gamma \rightarrow 1$ as $n$ grows
We numerically compute the values of the CBS cosntant $\gamma$, the condition number of the preconditioned operator and its estimate for a few values of $n$ in \cref{tab:example}.
We can see that as $n$ grows that $\gamma \rightarrow 1$ and that the condition number estimate is attained.

\begin{table}[htbp]
\centering
\begin{tabular}{l | l *{2}{|r}}
  $n$ & $\gamma$ & $\kappa_{eff} \left( \MM^{-1}\MA \right)$  & $\left( \kappa\left( \MA \right) +1  \right)^{2}/\left(4\kappa(\MA)\right)$\\\hline
  3 & $0.7071$ & $2.0000$ & $2.0000$\\\hline
  9 & $0.9511$ & $10.4721$ & $10.4721$\\\hline
  99 & $0.9995$ & $1,013.5452$ & $1,013.5452$\\\hline
  999 & $1.0000$ & $10,1321.5170$ & $10,1321.5170$
\end{tabular}
\caption{Numerically computed (and rounded) values of $\gamma$, the condition number of the preconditioned operator, and its estimate for the example with tightly coupled active and free sets.}
\label{tab:example}
\end{table}

We note that it is reasonable to expect that problems arising in many applications using finite element or similar discretizations will have rather weak coupling between the active and free sets for reasonable approximations of the solution due to the relative "smoothness" of the solution and constraints.

\subsection{Alternative Proof of Approximate Preconditioner Conditioning}
Here we present an alternative proof of Theorem \ref{theo:precond}, which does not use the CBS constant. However, it gives us the conditions on when the bound is sharp. The proof uses the Kantorovich inequality.

\begin{theorem}[Kantorovich inequality]
  \label{theo:kant}
  For any symmetric positive definite matrix $\MA \in \mathbb{R}^{n \times n}$, the following inequalities are valid:
\begin{equation*}
  1 \leq \frac{\left(\Vx^{T} \MA \Vx \right) \left(\Vx^{T} \MA^{-1} \Vx\right)}{\left(\Vx^{T} \Vx\right)^{2}} \leq \frac{\left( \kappa\left( \MA \right) +1  \right)^{2}}{4\kappa(\MA)} \quad \forall \Vx \in \mathbb{R}^{n}\setminus \left\{\Vo \right\}.
\end{equation*}
Moreover, the upper bound is sharp, i.e.,
\begin{equation*}
  \sup_{\Vx \neq \Vo} \frac{\left(\Vx^{T} \MA \Vx \right) \left(\Vx^{T} \MA^{-1} \Vx\right)}{\left(\Vx^{T} \Vx\right)^{2}} = \frac{\left( \kappa\left( \MA \right) +1  \right)^{2}}{4\kappa(\MA)}.
\end{equation*}
\end{theorem}
\begin{proof}
  See Lemma 3.11 of \cite{Axelsson96}.
\end{proof}
Interestingly, this is equivalent to the Wielandt inequality (\cref{theo:wielandt}) \cite{Zhang2001}, which we have used to bound the CBS constant $\gamma$ in the proof of the conditioning in \cref{theo:precond}.

\begin{theorem}[Conditioning of the preconditioned operator]
  \label{theo:condsimple}
  Let $\overline{\MM} = \MA$ and the inner preconditioner be the full inverse of $\overline{\MM}$. Then
    \begin{equation*}
      \kappa_{eff}\left(\MM^{-1}\MA\right) \leq \frac{\left( \kappa\left( \MA \right) +1  \right)^{2}}{4\kappa(\MA)}.
    \end{equation*}
\end{theorem}
\begin{proof}
  The maximal eigenvalue $\lambda_{max}$ of the preconditioned operator $\MM^{-1}\MA$ is by \cref{theo:eigvals} equal to the maximal eigenvalue of $\MS^{-1}\MA_{\mathcal{FF}}$, i.e.,
  \begin{equation*}
    \lambda_{max} = \sup_{\Vv_{\mathcal{F}} \neq \Vo} \frac{\Vv_{\mathcal{F}}^{T}\MA_{\mathcal{FF}}\Vv_{\mathcal{F}}}{\Vv_{\mathcal{F}}^{T}\MS_{\mathcal{FF}}\Vv_{\mathcal{F}}}.
  \end{equation*}
  From the Cauchy--Schwarz inequality
  \begin{equation}
    \label{eq:cs}
    \left(\Vv_{\mathcal{F}}^{T}\Vv_{\mathcal{F}}\right)^{2} = \left(\Vv_{\mathcal{F}}^{T}\MS^{\frac{1}{2}}\MS^{-\frac{1}{2}}\Vv_{\mathcal{F}}\right)^{2} \leq \left(\Vv_{\mathcal{F}}^{T}\MS\Vv_{\mathcal{F}}\right) \left(\Vv_{\mathcal{F}}^{T}\MS^{-1}\Vv_{\mathcal{F}}\right),
  \end{equation}
  we have
  \begin{equation*}
    \frac{1}{\Vv_{\mathcal{F}}^{T}\MS\Vv_{\mathcal{F}}} \leq \frac{\Vv_{\mathcal{F}}^{T}\MS^{-1}\Vv_{\mathcal{F}}}{\left(\Vv_{\mathcal{F}}^{T}\Vv_{\mathcal{F}}\right)^{2}}.
  \end{equation*}
  Multiplying both sides by $\Vv_{\mathcal{F}}^{T}\MA_{\mathcal{FF}}\Vv_{\mathcal{F}}$, we have for all $\Vv_{\mathcal{F}} \neq \Vo$
  \begin{equation*}
    \frac{\Vv_{\mathcal{F}}^{T}\MA_{\mathcal{FF}}\Vv_{\mathcal{F}}}{\Vv_{\mathcal{F}}^{T}\MS_{\mathcal{FF}}\Vv_{\mathcal{F}}} \leq
    \frac{\left(\Vv_{\mathcal{F}}^{T}\MA_{\mathcal{FF}}\Vv_{\mathcal{F}}\right) \left(\Vv_{\mathcal{F}}^{T}\MS^{-1}\Vv_{\mathcal{F}}\right)}{\left(\Vv_{\mathcal{F}}^{T}\Vv_{\mathcal{F}}\right)^{2}}.
  \end{equation*}
  Using $\left(\MA^{-1}\right)_\mathcal{FF} = \MS^{-1}$ from the block inversion formula and zero-padding $\Vv_{\mathcal{F}}$ as $\Vy = \begin{pmatrix} \Vv_{\mathcal{F}}\\ \Vo \end{pmatrix} \in \mathbb{R}^{n} \setminus\left\{\Vo\right\}$, we have
  \begin{equation*}
    \frac{\left(\Vv_{\mathcal{F}}^{T}\MA_{\mathcal{FF}}\Vv_{\mathcal{F}}\right) \left(\Vv_{\mathcal{F}}^{T}\MS^{-1}\Vv_{\mathcal{F}}\right)}{\left(\Vv_{\mathcal{F}}^{T}\Vv_{\mathcal{F}}\right)^{2}} =
    \frac{\left(\Vy^{T}\MA\Vy\right) \left(\Vy^{T}\MA^{-1}\Vy\right)}{\left(\Vy^{T}\Vy\right)^{2}}.
  \end{equation*}
  Since $\Vy$ is from a subset of $\mathbb{R}^{n}\setminus\left\{\Vo\right\}$, we bound the last term using the Kantorovich inequality (\cref{theo:kant}) by omitting the restriction on the zero components of $\Vy$ so that for all $\Vx \in \mathbb{R}^{n} \setminus\left\{\Vo\right\}$:
  \begin{equation*}
    \frac{\left(\Vy^{T}\MA\Vy\right) \left(\Vy^{T}\MA^{-1}\Vy\right)}{\left(\Vy^{T}\Vy\right)^{2}}
    \leq \frac{\left(\Vx^{T}\MA\Vx\right) \left(\Vx^{T}\MA^{-1}\Vx\right)}{\left(\Vx^{T}\Vx\right)^{2}}
    \leq \frac{\left( \kappa\left( \MA \right) +1 \right)^{2}}{4\kappa(\MA)}.
  \end{equation*}
  By \cref{theo:eigvals}, the minimal non-zero eigenvalue of the preconditioned operator $\MM^{-1}\MA$ is greater than or equal to $1$, which completes the proof.
\end{proof}

For the bound on the largest eigenvalue of the preconditioned operator to be sharp, both the Cauchy--Schwarz inequality \eqref{eq:cs} and the Kantorovich inequality must be sharp.
The bound in the Kantorovich inequality is attained for
\begin{equation*}
  \Vx = \Vv_{min} + \Vv_{max},
\end{equation*}
where $\Vv_{min}$ and $\Vv_{max}$ are eigenvectors belonging to the smallest and the largest eigenvalues of $\MA$, respectively.
The Cauchy--Schwarz inequality becomes equality if $\Vv_{\mathcal{F}}$ is an eigenvector of $\MS$.
Therefore, the final bound on the largest eigenvalue is sharp if $\Vv_{min}$ and $\Vv_{max}$ cancel on the active set, i.e.,
\begin{equation*}
\begin{pmatrix} \Vv_{\mathcal{F}}\\ \Vo \end{pmatrix} = \Vv_{min} + \Vv_{max},
\end{equation*}
and their sum on the free set is an eigenvector of $\MS$.

For example,
\begin{equation*}
  \Vx = \begin{pmatrix} 1 \\ 0 \end{pmatrix}, \qquad \MA = \begin{pmatrix} 2 & -1 \\ -1 & 2 \end{pmatrix},
\end{equation*}
will attain the supremum for both inequalities, but for
\begin{equation*}
\MA = \begin{pmatrix} 2 & -1 \\ -1 & 3 \end{pmatrix},
\end{equation*}
the supremum for the Kantorovich inequality is attained for
\begin{equation*}
  \Vx \approx \begin{pmatrix} -1.3764 \\ 0.3249 \end{pmatrix},
\end{equation*}
which cannot be zero on any non-trivial active set.

The bound was numerically attained in the experiments presented in \cref{tab:example} for the Laplace $n \times n$ tridiagonal matrix with $n \geq 3$ being and odd integer.
In fact, with the above discussion, we can show that the bound is attained analytically.

The $k$th component of the minimal and maximal eigenvectors are
\begin{equation*}
  \left( \Vv_{min} \right)_k = \sin\left( \frac{k\pi}{n+1} \right) \quad \text{and} \quad \left( \Vv_{max} \right)_k = \sin\left( \frac{nk\pi}{n+1} \right),
\end{equation*}
respectively. Rewriting $\Vv_{max}$ and using the angle difference identity we have
\begin{equation*}
  \left( \Vv_{max} \right)_k = \sin\left( \frac{nk\pi}{n+1} \right) = \sin\left( k\pi - \frac{k\pi}{n+1} \right)
  = \sin\left( k\pi \right) \cos\left( \frac{k\pi}{n+1} \right) - \cos\left( k\pi \right) \sin\left( \frac{k\pi}{n+1} \right).
\end{equation*}
Since $k$ is an integer, we have $\sin\left( k\pi \right) = 0$ and $\cos\left( k\pi \right) = \left( -1 \right) ^{k}$. Substituting above gives a simplified expression
\begin{equation*}
  \left( \Vv_{max} \right)_k = (-1)^{k+1} \sin\left( \frac{k\pi}{n+1} \right) = (-1)^{k+1} \left( \Vv_{min} \right)_k.
\end{equation*}
Therefore, with $\Vx = \Vv_{min} + \Vv_{max}$ the Kantorovich inequality is attained and $\Vx$ is precisely zero on our active set given by even integers.

It remains to show that $\Vx_{\mathcal{F}}$ is an eigenvector of the Schur complement $\MS$. The Schur complement for our problem is:
\begin{equation*}
  \MS = \MA_{\mathcal{FF}} - \MA_{\mathcal{FA}} \MA_{\mathcal{AA}}^{-1} \MA_{\mathcal{AF}} = 2\MI - \frac{1}{2}\MA_{\mathcal{FA}} \MA_{\mathcal{AF}}.
\end{equation*}
Because the coupling matrices $\MA_{\mathcal{FA}}$ and $\MA_{\mathcal{AF}}$ are upper and lower bidiagonal, respectively, with $-1$ entries, we have
\begin{equation*}
  \MS = 2\MI - \frac{1}{2}\begin{pmatrix} 1 & 1\\
  1 & 2 &1\\
    & \ddots & \ddots & \ddots\\
    & & 1 & 1
  \end{pmatrix}
  = \frac{1}{2} \begin{pmatrix} 3 & -1\\
  -1 & 2 &-1\\
    & \ddots & \ddots & \ddots\\
    & & -1 & 3
  \end{pmatrix}.
\end{equation*}
Since the $l$th component of $\Vx_\mathcal{F} \in \mathbb{R}^{\frac{n+1}{2}}$ is
\begin{equation*}
  \left( \Vx_\mathcal{F} \right)_l = 2 \sin \left( \frac{\left( 2l - 1 \right)\pi}{n+1} \right),
\end{equation*}
and letting
\begin{equation*}
  \varphi = \frac{\left( 2l - 1 \right)\pi}{n+1} \quad \text{and} \quad \delta = \frac{2\pi}{n+1}
\end{equation*}
we have for any row $l$ except the first and the last that
\begin{flalign*}
  \left(\MS \Vx_{\mathcal{F}}\right)_l
  &= -\sin\left( \varphi - \delta \right)
  + 2\sin \varphi
  - \sin\left( \varphi + \delta \right)
  = 2\sin \varphi
  - \left[ \sin\left( \varphi - \delta \right) + \sin\left( \varphi +\delta \right)\right]\\
  &= 2\sin \varphi - \left[ 2\sin \varphi \cos\delta \right] = 2\left(1- \cos\delta  \right) \sin \varphi = \left(1- \cos\delta  \right)\left(\Vx_{\mathcal{F}}\right)_l,
\end{flalign*}
where we have used the angle sum and difference identities and that cosine is even.
For the first row $l=1$ we have
\begin{flalign*}
  \left(\MS \Vx_{\mathcal{F}}\right)_1
  &= 3 \sin \varphi - \sin 3\varphi
  = 3 \sin \varphi - \left[ 3\sin\varphi - 4\sin^{3}\varphi \right]
  = 2\left( 1 -\cos 2\varphi \right) \sin \varphi\\
  &= 2\left(1- \cos\delta  \right) \sin \varphi =
  \left(1- \cos\delta  \right)\left(\Vx_{\mathcal{F}}\right)_1,
\end{flalign*}
where we have used the triple-angle and double-angle formulas.
Finally, for the last row $l = \frac{n+1}{2}$ we have
\begin{flalign*}
  \left(\Vx_{\mathcal{F}}\right)_l
  &= \sin \frac{n\pi}{n+1} = \sin\left( \pi - \frac{\pi}{n+1} \right) = \sin \frac{\pi}{n+1} = \left(\Vx_{\mathcal{F}}\right)_1\\
  \left(\Vx_{\mathcal{F}}\right)_{l-1},
  &= \sin \frac{\left( n-2 \right) \pi}{n+1} = \sin\left( \pi - 3\frac{\pi}{n+1} \right)  = \sin 3\frac{\pi}{n+1} = \left(\Vx_{\mathcal{F}}\right)_2,
\end{flalign*}
so that
\begin{equation*}
\left(\MS \Vx_{\mathcal{F}}\right)_l = \left(\MS \Vx_{\mathcal{F}}\right)_1 = \left(1- \cos\delta  \right)\left(\Vx_{\mathcal{F}}\right)_1 = \left(1- \cos\delta \right)\left(\Vx_{\mathcal{F}}\right)_l.
\end{equation*}
Therefore, $\Vx_{\mathcal{F}}$ is an eigenvector of the Schur complement $\MS$ and the Cauchy--Schwarz inequality bound is attained.
Since $\begin{pmatrix} 1 & -1 & \dots \end{pmatrix}$ is in the null space of $\MA_{\mathcal{AF}}$ the minimal non-zero eigenvalue of the preconditioned operator is $1$ by \cref{theo:eigvals}, and because both the Cauchy--Schwarz and Kantorovich bounds are attained, the bound on the preconditioned operator is attained.

\section{Related Work}
\label{sec:related}
As far as we know, these are the only results of the MPGP-type algorithm showcasing the preconditioning in face (results for partially constrained problems using deflation can be found in \cite{Domoradova2006,Dos-book-09,Jarosova2012}).

The idea of the approximate preconditioning used for MPRGP can first be found in the accompanying codes to the article by Narain et al.\cite{Narain2010}, where it was used in combination with the incomplete Cholesky preconditioner.
The article does not contain any details about the MPRGP algorithm, the preconditioning in face, nor any numerical experiments related to MPRGP and its preconditioning.
A related work by Gerszewski and Bargteil \cite{Gerszewski2013} uses MPRGP preconditioned by the incomplete Cholesky.
While the article cites Narain et al., it is not obvious which variant of the preconditioning in face is used.
In any case, there is again no research presented with respect to the preconditioned method.

Finally, a variant of the approximate preconditioning in face paired with MPPCG is used by Takahashi and Batty in \cite{Takahashi2023}.
In our notation, they assemble aggregation-based algebraic multigrid for the full Hessian as the inner preconditioner, but instead of restricting the preconditioning only to the free set, they filter all indices that are connected to the active components through the restriction operators.
Indeed, restricting inner preconditioner application to indices that are not connected by the inner preconditioner application to any active components should intuitively provide better efficiency of the preconditioning - see the discussion of the CBS constant based condition number estimate after \cref{theo:precond}.

\FloatBarrier
\section{Numerical Experiments}
\label{sec:results}
The open-source PERMON library \cite{permon,permongit} was used to compute the numerical experiments.
PERMON stands for Parallel, Efficient, Robust, Modular, Object-oriented, Numerical.
It provides solvers and a number of transformations and other helpful functions for the solution of QP problems, as well as FETI-type domain decomposition methods and support vector machines.
PERMON is built on top of PETSc \cite{petsc-web-page,petsc-user-ref,petsc-efficient}, utilizing the same programming style. Therefore, it is written in C, uses MPI for parallelization, and is able to utilize GPUs/accelerators from a growing number of vendors.

Standard preconditioners available in PETSc with the default options were used to compute the results.
The Cholesky "preconditioner" represents the application of the direct solver, i.e., preconditioning by the inverse of the Hessian, using MUMPS \cite{Amestoy2001,Amestoy2019}.
ICC is the incomplete Cholesky factorization \cite{Chan1997}, and SSOR is the symmetric successive over-relaxation \cite{Young1971}.

The CG method, applied to the system of linear equations preconditioned by the inverse of the Hessian, will converge in a single iteration.
That is not the case for the preconditioned MPGP-type methods because the active set needs to be identified.
However, if we start with the correct active set or once the correct active set is identified, the preconditioned method converges to machine precision in a single iteration.
In any case, since the inverse preconditioner is the optimal preconditioner for the preconditioning in face in the sense of preconditioning the Hessian on the free set, the number of Hessian multiplications for the Cholesky preconditioner is a very interesting metric\footnote{
Although it might not be the lower bound on the number of Hessian multiplications needed by the preconditioned MPGP-type algorithm, despite this being the case in our numerical experiments.}.

In the unpreconditioned method, the multiplication by the Hessian will typically take over 90 \% of the runtime.
Additionally, the preconditioning affects the number of CG, expansion, and proportioning steps, and the expansion step needs two Hessian multiplications, while the other two steps need only one.
Therefore, the number of Hessian multiplications is a better metric than the number of iterations to assess the numerical behavior of MPGP-type methods.
In the results presented below, we report the number of Hessian multiplications as well as the number of each step.
The number of iterations can be computed as the sum of the number of each step.

Due to the inclusion of the preconditioner in each iteration, the number of Hessian multiplications, while still of interest, cannot be used as a metric for comparison between the methods.
Therefore, timings and speedups based on the timings are provided.
To ensure a high quality of timings, the presented results were computed with an optimized build\footnote{The code and libraries are built by Cray clang 16 with -O3 flag} on a single core of a dedicated node of the LUMI supercomputer \cite{lumi}, i.e., on AMD EPYC 7763 at 2.45 GHz.
The stopping criterion was the relative tolerance of $10^{-10}$.

The first problem is a variant of determining the pressure distribution of the journal bearing problem from the MINPACK-2 test problem collection \cite{Averick1992}.
This 2D problem corresponds to tutorial \emph{jbearing2} in PERMON.
The second problem is a 3D linear elastic cube that is fixed at the bottom, pushed from the top, and there is an obstacle parallel to the right face at a small distance away, which results in an upper bound constraint on displacement.
The problems are discretized by the P1 and Q1 Lagrange finite elements, respectively.
Complete descriptions of the problems with all parameters are available in \cite{KruzikPhD}.

A number of increasingly refined discretizations is presented for each problem as this gives an interesting comparison as the conditioning of the Hessian deteriorates.
The results for the 3D linear elasticity are in \cref{tab:precondCube,tab:precondCube2} and for the journal bearing problem in \cref{tab:precondJb,tab:precondJb2,tab:precondJb3,tab:precondJb4}.
Columns S$_b$ and S$_M$ contain speedups.
S$_b$ is computed with respect to the same unpreconditioned method, while S$_M$ is computed with respect to the unpreconditioned MPRGP method.

First, examining the performance of the standard MPRGP with the preconditioning in face, the number of Hessian multiplications is significantly reduced compared to the unpreconditioned method.
This reduction is driven by a large decrease in the CG steps, which is precisely what we would expect.
The number of expansion steps appears to be proportional to the preconditioner effectiveness, being the lowest for the Cholesky preconditioner, followed by ICC, and finally SSOR.
Compared to the unpreconditioned method, the number of expansions was typically lower for the journal bearing problem and higher for the elasticity problem.
The number of proportioning steps follows similar trends as the expansion steps, but the change between the preconditioners is much less pronounced.

As for the new approximate preconditioning in face combined with the standard MPRGP, the effectiveness in reducing the number of CG iterations is still there.
However, there is a further increase in the number of expansion steps, which is very noticeable in journal bearing problems.
The cause of the increase could be driven by the decrease in the effectiveness of the approximate preconditioning in face compared to the standard preconditioning in face.
Overall, the number of Hessian multiplications usually increases compared to the standard preconditioning in face, especially for larger journal bearing problems.
Despite this, the time needed by the approximate preconditioning in face is significantly lower than the preconditioning in face (except for the ICC preconditioner in \cref{tab:precondJb4}, where the preconditioning in face is slightly faster).
The variants of the preconditioners in face applying the ICC preconditioner exhibited a speedup between $1.96$ and $4.66$ for the preconditioning in face, and between $4.28$ and $6.43$ for the approximate variant on the journal bearing problem.
However, they were much slower on the elasticity benchmark, where they attained a speedup of at most $0.15$ and $0.86$ for the preconditioning in face and its approximate variant, respectively.

Despite the preconditioning working, i.e., the number of CG steps is reduced, the growth in the number of expansion steps limits the usefulness of the preconditioning.
Fortunately, we have the MPPCG method that was specifically designed to reduce the number of expansions.
The results show that the MPPCG method significantly limits the number of expansion steps, while the number of CG and proportioning steps is in the same ballpark, if not nearly identical, compared to the MPRGP method.
Even the unpreconditioned MPPCG method exhibits some speedup over the unpreconditioned MPRGP.
The preconditioning in face always performs better than the approximate variant in terms of the number of Hessian multiplications, but worse in terms of the time to solution.
The approximate preconditioning in face exhibits small speedups even for the SSOR preconditioner, ranging from $1.14$ to $1.94$.
Equipping the approximate preconditioner with ICC leads to very large speedups between $2.70$ and $10.38$.
If the unpreconditioned MPRGP is taken as the base, then the speedups range between $5.13$ and $13.46$.

\begin{table}[htbp]
\centering
\begin{tabular}{l *{2}{|l} *{7}{|r}}
  Method & Type & Precond. & Hess. & CG & Exp. & Prop. & Time [s] & $S_b$ & $S_M$\\\hline
  MPRGP & None & None & 2030 & 1262 & 381 & 5 & 1.81 & 1.00 & 1.00\\ \hline
MPRGP & Face & Cholesky & 788 & 5 & 389 & 4 & 227.16 & 0.01 & 0.01\\
MPRGP & Approx & Cholesky & 817 & 28 & 392 & 4 & 9.51 & 0.19 & 0.19\\ \hline
MPRGP & Face & ICC & 1154 & 95 & 526 & 6 & 20.65 & 0.09 & 0.09\\
MPRGP & Approx & ICC & 1357 & 99 & 626 & 5 & 2.17 & 0.84 & 0.84\\ \hline
MPRGP & Face & SSOR & 1617 & 178 & 717 & 4 & 15.94 & 0.11 & 0.11\\
MPRGP & Approx & SSOR & 1642 & 191 & 723 & 4 & 2.31 & 0.78 & 0.78\\ \hline \hline
MPPCG & None & None & 1054 & 864 & 92 & 5 & 0.95 & 1.00 & 1.90\\ \hline
MPPCG & Face & Cholesky & 12 & 5 & 1 & 4 & 6.24 & 0.15 & 0.29\\
MPPCG & Approx & Cholesky & 35 & 28 & 1 & 4 & 1.22 & 0.78 & 1.48\\ \hline
MPPCG & Face & ICC & 117 & 83 & 14 & 5 & 3.26 & 0.29 & 0.56\\
MPPCG & Approx & ICC & 163 & 127 & 15 & 5 & 0.35 & 2.70 & 5.13\\ \hline
MPPCG & Face & SSOR & 257 & 192 & 30 & 4 & 3.73 & 0.26 & 0.49\\
MPPCG & Approx & SSOR & 285 & 202 & 39 & 4 & 0.49 & 1.94 & 3.68\\
\end{tabular}
\caption[Preconditioned 3D contact problem I.]{Results for preconditioning the 3D linear elasticity cube contact problem with 10x20x40 finite elements (28,413 DOFs).}
\label{tab:precondCube}
\end{table}

\begin{table}[htbp]
\centering
\begin{tabular}{l *{2}{|l} *{7}{|r}}
  Method & Type & Precond. & Hess. & CG & Exp. & Prop. & Time [s] & $S_b$ & $S_M$\\\hline
MPRGP & None & None & 6544 & 3590 & 1472 & 9 & 88.51 & 1.00 & 1.00\\ \hline
MPRGP & Face & Cholesky & 1818 & 6 & 903 & 5 & 14883.00 & 0.01 & 0.01\\
MPRGP & Approx & Cholesky & 3095 & 44 & 1522 & 6 & 439.77 & 0.20 & 0.20\\ \hline
MPRGP & Face & ICC & 4258 & 209 & 2020 & 8 & 594.58 & 0.15 & 0.15\\
MPRGP & Approx & ICC & 5446 & 350 & 2544 & 7 & 102.93 & 0.86 & 0.86\\ \hline
MPRGP & Face & SSOR & 5964 & 371 & 2793 & 6 & 487.90 & 0.18 & 0.18\\
MPRGP & Approx & SSOR & 6040 & 405 & 2814 & 6 & 152.52 & 0.58 & 0.58\\ \hline \hline

MPPCG & None & None & 2766 & 2269 & 244 & 8 & 37.93 & 1.00 & 2.33\\ \hline
MPPCG & Face & Cholesky & 14 & 6 & 1 & 5 & 212.37 & 0.18 & 0.42\\
MPPCG & Approx & Cholesky & 57 & 43 & 3 & 7 & 30.19 & 1.26 & 2.93\\ \hline
MPPCG & Face & ICC & 344 & 212 & 60 & 11 & 72.38 & 0.52 & 1.22\\
MPPCG & Approx & ICC & 473 & 297 & 84 & 7 & 10.38 & 3.65 & 8.53\\ \hline
MPPCG & Face & SSOR & 696 & 439 & 125 & 6 & 82.46 & 0.46 & 1.07\\
MPPCG & Approx & SSOR & 715 & 443 & 132 & 7 & 22.55 & 1.68 & 3.93\\
\end{tabular}
\caption[Preconditioned 3D contact problem II.]{Results for preconditioning the 3D linear elasticity cube contact problem with 20x40x80 finite elements (209,223 DOFs).}
\label{tab:precondCube2}
\end{table}

\begin{table}[htbp]
\centering
\begin{tabular}{l *{2}{|l} *{7}{|r}}
  Method & Type & Precond. & Hess. & CG & Exp. & Prop. & Time [s] & $S_b$ & $S_M$\\\hline
MPRGP & None & None & 2884 & 2660 & 69 & 85 & 0.33 & 1.00 & 1.00\\ \hline
MPRGP & Face & Cholesky & 157 & 78 & 0 & 78 & 1.95 & 0.17 & 0.17\\
MPRGP & Approx & Cholesky & 494 & 79 & 165 & 84 & 1.63 & 0.20 & 0.20\\ \hline
MPRGP & Face & ICC & 179 & 100 & 0 & 78 & 0.17 & 1.96 & 1.96\\
MPRGP & Approx & ICC & 308 & 79 & 72 & 84 & 0.05 & 6.43 & 6.43\\ \hline
MPRGP & Face & SSOR & 999 & 732 & 93 & 80 & 0.77 & 0.43 & 0.43\\
MPRGP & Approx & SSOR & 994 & 666 & 122 & 83 & 0.21 & 1.59 & 1.59\\ \hline \hline

MPPCG & None & None & 2348 & 2218 & 25 & 79 & 0.27 & 1.00 & 1.24\\ \hline
MPPCG & Face & Cholesky & 157 & 78 & 0 & 78 & 1.95 & 0.14 & 0.17\\
MPPCG & Approx & Cholesky & 197 & 97 & 10 & 79 & 0.91 & 0.29 & 0.36\\ \hline
MPPCG & Face & ICC & 179 & 100 & 0 & 78 & 0.17 & 1.59 & 1.96\\
MPPCG & Approx & ICC & 208 & 87 & 19 & 82 & 0.04 & 7.28 & 9.01\\ \hline
MPPCG & Face & SSOR & 748 & 623 & 22 & 80 & 0.60 & 0.44 & 0.55\\
MPPCG & Approx & SSOR & 858 & 699 & 38 & 82 & 0.19 & 1.42 & 1.76\\
\end{tabular}
\caption[Preconditioned journal bearing problem I.]{Results for preconditioning the journal bearing problem with 400x25 discretization points (10,000 DOFs).}
\label{tab:precondJb}
\end{table}

\begin{table}[htbp]
\centering
\begin{tabular}{l *{2}{|l} *{7}{|r}}
  Method & Type & Precond. & Hess. & CG & Exp. & Prop. & Time [s] & $S_b$ & $S_M$\\\hline
MPRGP & None & None & 7789 & 6989 & 306 & 187 & 3.30 & 1.00 & 1.00\\ \hline
MPRGP & Face & Cholesky & 309 & 154 & 0 & 154 & 35.22 & 0.09 & 0.09\\
MPRGP & Approx & Cholesky & 856 & 150 & 274 & 157 & 10.50 & 0.31 & 0.31\\ \hline
MPRGP & Face & ICC & 366 & 189 & 11 & 154 & 1.29 & 2.56 & 2.56\\
MPRGP & Approx & ICC & 1092 & 128 & 379 & 205 & 0.65 & 5.11 & 5.11\\ \hline
MPRGP & Face & SSOR & 3168 & 1633 & 667 & 200 & 8.29 & 0.40 & 0.40\\
MPRGP & Approx & SSOR & 4928 & 2344 & 1173 & 237 & 3.71 & 0.89 & 0.89\\ \hline \hline
MPPCG & None & None & 7286 & 6578 & 260 & 187 & 3.03 & 1.00 & 1.09\\ \hline
MPPCG & Face & Cholesky & 309 & 154 & 0 & 154 & 35.16 & 0.09 & 0.09\\
MPPCG & Approx & Cholesky & 421 & 196 & 33 & 158 & 7.09 & 0.43 & 0.47\\ \hline
MPPCG & Face & ICC & 352 & 191 & 3 & 154 & 1.26 & 2.40 & 2.61\\
MPPCG & Approx & ICC & 454 & 154 & 64 & 171 & 0.29 & 10.38 & 11.30\\ \hline
MPPCG & Face & SSOR & 2066 & 1538 & 159 & 209 & 6.14 & 0.49 & 0.54\\
MPPCG & Approx & SSOR & 2823 & 1970 & 308 & 236 & 2.25 & 1.35 & 1.47\\
\end{tabular}
\caption[Preconditioned journal bearing problem II.]{Results for preconditioning the journal bearing problem with 800x50 discretization points (40,000 DOFs).}
\label{tab:precondJb2}
\end{table}

\begin{table}[htbp]
\centering
\begin{tabular}{l *{2}{|l} *{7}{|r}}
  Method & Type & Precond. & Hess. & CG & Exp. & Prop. & Time [s] & $S_b$ & $S_M$\\\hline
MPRGP & None & None & 12022 & 9389 & 1199 & 234 & 9.95 & 1.00 & 1.00\\ \hline
MPRGP & Face & Cholesky & 309 & 154 & 0 & 154 & 75.85 & 0.13 & 0.13\\
MPRGP & Approx & Cholesky & 1839 & 148 & 765 & 160 & 39.96 & 0.25 & 0.25\\ \hline
MPRGP & Face & ICC & 507 & 222 & 64 & 156 & 3.32 & 3.00 & 3.00\\
MPRGP & Approx & ICC & 1920 & 140 & 772 & 235 & 2.16 & 4.60 & 4.60\\ \hline
MPRGP & Face & SSOR & 4634 & 1971 & 1226 & 210 & 22.85 & 0.44 & 0.44\\
MPRGP & Approx & SSOR & 7330 & 2667 & 2185 & 292 & 10.45 & 0.95 & 0.95\\ \hline \hline
MPPCG & None & None & 8906 & 7809 & 440 & 216 & 7.39 & 1.00 & 1.35\\ \hline
MPPCG & Face & Cholesky & 309 & 154 & 0 & 154 & 75.86 & 0.10 & 0.13\\
MPPCG & Approx & Cholesky & 459 & 208 & 46 & 158 & 15.35 & 0.48 & 0.65\\ \hline
MPPCG & Face & ICC & 457 & 217 & 38 & 163 & 3.10 & 2.38 & 3.21\\
MPPCG & Approx & ICC & 1042 & 169 & 274 & 324 & 1.17 & 6.29 & 8.47\\ \hline
MPPCG & Face & SSOR & 2853 & 1913 & 357 & 225 & 16.24 & 0.46 & 0.61\\
MPPCG & Approx & SSOR & 2996 & 1961 & 409 & 216 & 4.64 & 1.59 & 2.14\\
\end{tabular}
\caption[Preconditioned journal bearing problem III.]{Results for preconditioning the journal bearing problem with 800x100 discretization points (80,000 DOFs).}
\label{tab:precondJb3}
\end{table}

\begin{table}[htbp]
\centering
\begin{tabular}{l *{2}{|l} *{7}{|r}}
  Method & Type & Precond. & Hess. & CG & Exp. & Prop. & Time [s] & $S_b$ & $S_M$\\\hline
MPRGP & None & None & 37044 & 28703 & 3844 & 652 & 60.14 & 1.00 & 1.00\\ \hline
MPRGP & Face & Cholesky & 617 & 308 & 0 & 308 & 317.32 & 0.19 & 0.19\\
MPRGP & Approx & Cholesky & 3612 & 244 & 1525 & 317 & 156.26 & 0.38 & 0.38\\ \hline
MPRGP & Face & ICC & 987 & 357 & 159 & 311 & 12.91 & 4.66 & 4.66\\
MPRGP & Approx & ICC & 6225 & 250 & 2738 & 498 & 14.06 & 4.28 & 4.28\\ \hline
MPRGP & Face & SSOR & 14072 & 5986 & 3780 & 525 & 144.25 & 0.42 & 0.42\\
MPRGP & Approx & SSOR & 25871 & 8442 & 8281 & 866 & 73.02 & 0.82 & 0.82\\ \hline \hline
MPPCG & None & None & 25166 & 21632 & 1509 & 515 & 40.40 & 1.00 & 1.49\\ \hline
MPPCG & Face & Cholesky & 617 & 308 & 0 & 308 & 317.43 & 0.13 & 0.19\\
MPPCG & Approx & Cholesky & 887 & 379 & 93 & 321 & 59.38 & 0.68 & 1.01\\ \hline
MPPCG & Face & ICC & 776 & 368 & 42 & 323 & 11.15 & 3.62 & 5.40\\
MPPCG & Approx & ICC & 1976 & 238 & 564 & 609 & 4.47 & 9.04 & 13.46\\ \hline
MPPCG & Face & SSOR & 9609 & 6194 & 1346 & 722 & 113.18 & 0.36 & 0.53\\
MPPCG & Approx & SSOR & 11661 & 6902 & 1982 & 794 & 35.32 & 1.14 & 1.70\\
\end{tabular}
\caption[Preconditioned journal bearing problem IV.]{Results for preconditioning the journal bearing problem with 1600x100 discretization points (160,000 DOFs).}
\label{tab:precondJb4}
\end{table}

\section{Conclusion}
%TODO improve conclusion
Approximate preconditioning in face for MPGP-type algorithms has been presented.
The main advantage of the approximate preconditioning in face is that the inner preconditioner needs to be computed only once, as opposed to on every change of the active/free sets in the case of the standard preconditioning in face.
This results in the approximate variant being much cheaper but typically requiring more Hessian multiplications, which are somewhat equivalent to the number of iterations in other algorithms.
The difference with the standard preconditioning in face has been demonstrated both numerically and, in specific cases, analytically.
We have also provided a sharp bound on the condition number of the preconditioned operator and constructed an example demonstrating the worst case bound.

The numerical experiments suggest that the approximate preconditioning in face suffers from an increase in the number of expansion steps.
This increase can be significantly reduced by the use of the MPPCG variant, which uses the projected conjugate gradient step for the expansion of the active set.
Overall, the observed speedup between the unpreconditioned MPPCG and MPPCG with preconditioning in face applying the best inner preconditioner ranges between $0.29$ to $3.62$.
On the other hand, the approximate preconditioning in face gives speedups between $2.70$ and $10.38$. When compared with the unpreconditioned MPRGP, the MPPCG method with the approximate preconditioning in face gives even larger speedups between $5.13$ and $13.46$.

In the future, we would like to apply the MPPCG method equipped with the approximate preconditioning in face to the solution of QP problems with known good or even optimal preconditioners for the unconstrained problem.
Such problems are, for example, contact problems in mechanics solved using the FETI method.

\section*{Acknowledgements}
The authors acknowledge the financial support of the European Union under the REFRESH - Research Excellence For Region Sustainability and High-tech Industries project number CZ.10.03.01/00/22 003/0000048 via the Operational Programme Just Transition.
This work was also supported by the European Union through the Operational Programme Jan Amos Komenský under project INODIN, number CZ.02.01.01/00/23 020/0008487.

\printbibliography
\end{document}